\documentclass[12pt]{article}
\usepackage{amssymb,latexsym, amsmath, amscd, array, graphicx}
\usepackage{amssymb}

\usepackage{setspace}

\usepackage{circledsteps}

\input xy
\xyoption{all}
\makeatletter 

\usepackage{setspace}
\usepackage{enumitem}

\date{}

\newtheorem{theorem}{Theorem}[section]

\newtheorem{proposition}[theorem]{Proposition}

\newcommand{\z}{{\Bbb Z}}

\newcommand{\re}{{\Bbb R}}
\newcommand{\N}{{\Bbb N}}

\newcommand{\lo}{\rightarrow}

\newcommand{\black}{{\blacksquare}}

\newcommand{\diam}{{\rm diam}}

\newcommand{\mdim}{{\rm mdim}}
\newcommand{\mesh}{{\rm mesh}}
\newcommand{\ord}{{\rm ord}}
 \newcommand{\dg}{{\dagger}}

\begin{document}

\title{\bf A union theorem  for mean dimension}

\author{  Michael Levin
\footnote{The author
 was supported by the Israel Science Foundation grant No. 2196/20 and 
 the  Fields institute of  mathematics  where 
 this note  was written  during his visit in 2023.}}

\maketitle

\begin{abstract}{Let $(X,\z)$ be a dynamical system on a compact metric $X$ and let
$X=\cup_{i\in \N} X_i$ be the countable union of closed invariant subsets $X_i$. We prove that
$\mdim X =\sup \{ \mdim X_i : i\in \N\}$.}
\\
{\bf Math. Subj. Class.:}  37B05 (54F45)
\end{abstract}
\begin{section}{Introduction}
The main goal of this note is to prove
\begin{theorem}
\label{theorem-union}
Let $(X,\z)$ be a dynamical system on a compact metric $X$ and let
$X=\cup_{i\in \N} X_i$ be the countable union of closed invariant subsets $X_i$. Then
$\mdim X =\sup \{ \mdim X_i : i\in \N\}$.
\end{theorem}
This theorem  is suggested by a similar result  for covering dimension and
it is quite  surprising  and was brought to the author's attention by Tom Meyerovitch  that  seemingly unrelated techniques
introduced in \cite{levin-1} can be used in proving Theorem \ref{theorem-union}.
  We adopt the notations of \cite{levin-1} and  the present note is based
on the   following results   of \cite{levin-1} and   ideas  used in their proofs.

\begin{theorem}
\rm{ (\cite{levin-1})}
\label{theorem-borel}
For any dynamical system $(X,\z)$ on a compact metric $X$  one has  \\
$\mdim X \times_\z \re=\mdim X$ where $X \times_\z \re$  is  Borel's construction for $(X,\z)$.
\end{theorem}

\begin{theorem}
\rm{ (\cite{levin-1})}
\label{theorem-minimal}
 Let $(X, \z)$ be a minimal dynamical system on a compact metric $X$ with
 $\mdim X =d$  and let $k>d$ be a natural number. Then   almost every map $f : X \lo [0,1]^k$
 induces the map $f^\z : X \lo  ([0,1]^k)^\z$ whose fibers  contain
 at most $[\frac{k}{k-d}]\frac{k}{k-d}$ points. 

\end{theorem}
Although  the proof of  Theorem \ref{theorem-union}  relies on Theorem \ref{theorem-borel}, 
we point out that Theorem \ref{theorem-union} trivially implies Theorem \ref{theorem-borel}
because one can easily split  $Y=X \times_\z \re$ into two invariant closed subsets $Y=Y_1 \cup Y_2$
with $\mdim Y_1 = \mdim Y_2 =\mdim X$.

Theorem \ref{theorem-minimal} falls short of our needs. One can show that the requirement 
that $(X,\z)$ in Theorem \ref{theorem-minimal} is minimal can be  weakened to
having the marker property that is enough for proving Theorem \ref{theorem-union}. 
 However such an extension of Theorem \ref{theorem-union} requires  a more subtle and complicated argument
 and will be presented elsewhere. In this note we prove the following  simplified version that meets
 our needs.

 \begin{theorem}
 \label{theorem-extension}
 Let $\mu: (X,\z) \lo (M, \z)$ be a surjective  equivariant  map of dynamical systems on compact metric spaces such that
 $(M, \z)$ is minimal and non-trivial,  let $d=\mdim (X,\z)$ and  let  $k>d$ be a natural number. Then for almost every map $f : X \lo [0,1]^k$
 the fibers of  the map
 $(\mu, f^\z) : X \lo  M \times   ([0,1]^k)^\z $  contain at most $[\frac{k}{k-d}]$ points. 
 \end{theorem}

It was conjectured in \cite{levin-1} and still remains open 
  that   the right estimate for the size of the fibers of $f^\z$ in Theorem \ref{theorem-minimal}
should   be  $[\frac{k}{k-d}]$.  Theorem \ref{theorem-extension} seems to be an indication 
supporting this conjecture. Moreover, since any non-trivial minimal system is an extension of a non-trivial minimal system
of an arbitrarily small mean dimension we may   assume in Theorem \ref{theorem-extension} 
 that $\mdim M < 1/2$ and hence  $M$ is embeddable in $[0,1]^\z$ and  get
 that $X$   admits an equivariant map to $([0,1]^{k+1})^\z$ 
whose fibers contain at most $[\frac{k}{k-d}]$ points, the estimate  in  some cases 
 better than the one provided by Theorem  \ref{theorem-minimal}.
 
 Let us finally mention that the requirement  in Theorem \ref{theorem-extension}  that $M$ is minimal and non-trivial can
 be weakened to $M$ having the marker property and the proof of Theorem \ref{theorem-extension} works for this case
 without any change.

 In the next section we show how to derive Theorem \ref{theorem-union} from Theorem \ref{theorem-minimal}
 and in the last section we prove Theorem \ref{theorem-minimal}.

\end{section}
\begin{section}{Proof of Theorem \ref{theorem-union}}
We show  here  how to derive Theorem \ref{theorem-union} from  Theorem \ref{theorem-extension} and will need the following
 facts. Recall that a map between topological spaces is said to be $0$-dimensional
if the fibers of the map are  of covering dimension $0$.
\begin{theorem}
\label{theorem-light-map}
Let  $f : (X,\z) \lo (Y,\z)$ be a $0$-dimensional equivariant map of dynamical system 
on compact metric spaces. Then $\mdim X \leq \mdim Y$.
\end{theorem}
{\bf Proof.} Let $\epsilon>0$. Since $f$ is $0$-dimensional there is $\delta >0$ 
 such that for every closed subset $A \subset Y$  of $\diam < \delta$
every connected component of $f^{-1}(A)$  is of $\diam < \epsilon$.  Clearly the theorem holds for $\mdim Y=\infty$.
Assume that $\mdim Y <\infty$ and let $d>0$ be  a real number such that $\mdim Y < d$. Take a finite closed
cover $\cal  A$ of $Y$ such that $\mesh ({\cal A}+z) < \delta$ for every integer $0\leq zd < \ord \cal A$.
Then for every $A \in \cal A$ and every connected component $C$  of $f^{-1}(A)$ we have that
$\diam (C +z) < \epsilon$ for every integer $ 0\leq zd <\ord \cal A$.  This implies that for each $A \in \cal A$
the set  $f^{-1}(A)$ splits into a finite family ${\cal B}_A$ of closed disjoint subsets of $X$
such that $\mesh ({\cal B}_A +z) <\epsilon $  for every integer $0\leq zd < \ord \cal A$.
Set $\cal B$ to be the union of ${\cal B}_A$ for all $A \in \cal A$. Then  $\cal B$ covers $X$ and $\ord {\cal B} \leq  \ord  {\cal A}$.
Thus we have  $\mesh ({\cal B}+z) < \epsilon$ for every integer $0\leq zd < \ord \cal B$ and hence 
$\cal B$ witnesses that $\mdim X \leq  d$ and the theorem follows. $\black$

\begin{proposition}
\label{proposition-selection}
Let $f : (X,\z)\lo (Y,\z)$ be an equivariant map of dynamical systems on compact metric spaces that admits
 a continuous (not necessarily equivariant) selection, i.e. there is a map  $s : Y \lo X $ such that  $f(s(y))=y $ for every $y \in Y$.
 Then $\mdim X \geq \mdim Y$.
 \end{proposition}
 {\bf Proof.} Let $\mdim X < d$. Take any  $\epsilon >0$ and let $\delta >0$ be such that the image under $f$ of any subset of
  $\diam < \delta$ in $X$  is of $\diam < \epsilon$ in $Y$. Since $\mdim  X < d$ there is a closed finite  cover ${\cal F}_X$ of $X$ such that
  $\mesh ({\cal F}_X +z) <\delta$ for every integer $z$ satisfying $0\leq zd  < \ord {\cal F}_X$. Let ${\cal S}$ be the cover ${\cal F}_X$ 
  restricted to $s(Y)$ and let ${\cal F}_Y=f({\cal S})$. Then $\ord {\cal F}_Y=\ord {\cal S} \leq \ord {\cal F}_X$ and
  hence $\mesh ({\cal S}+z) \leq \mesh ({\cal F}_X +z)<\delta$ 
  for  every integer $z$ satisfying  $0\leq zd\leq \ord {\cal S}$. Since ${\cal F}_Y +z=f({\cal S}+z)$
  we have $\mesh({\cal F}_Y +z) < \epsilon$ for every integer $z$ satisfying  $0\leq zd< \ord {\cal F}_Y$. This implies $\mdim Y < d$
  and the proposition follows. $\black$
\\\\
{\bf Proof of Theorem \ref{theorem-union}.} Take any non-trivial minimal  dynamical  system $(M, \z)$ on a metric compact  $M$ with
$\mdim M =0$. Let $n$ be any positive integer.
Consider the dynamical system $(X,n\z)$ with the action of the subgroup $n\z$ of $\z$ and 
consider the dynamical system $(Y,\z)=(X,n\z)\times (M,\z)$ with the product action and
let   $\mu$  be
the projection of $(Y,\z)$ to $(M,\z)$.  
 By Proposition \ref{proposition-selection}
   $\mdim (X, n\z)  \leq \mdim (Y,\z)$. 
   Then, since  $\mdim (Y,\z) \leq \mdim (X,n\z) +\mdim (M,\z)=\mdim (X, n\z)=n\mdim(X, \z)$ \cite{lindenstrauss-weiss},
    we get   $\mdim (Y,\z)=  n\mdim (X, \z)$.
Similarly we conclude that $\mdim  (Y_i, \z) =n\mdim (X_i, \z)$ for every $i$ where $(Y_i, \z) = (X_i, n\z)\times (M, \z)$. 

Let $d>0$ be such that  $\mdim (X_i, \z) < d$ for every $i$ and  let $k=[nd]+1$.
By Theorem \ref{theorem-extension}
there is a map $f : Y \lo [0,1]^k$  which induces  the map $f^\z : (Y,\z) \lo ([0,1]^k)^\z$ such that the fibers of
$(\mu, f^\z) : (Y, \z) \lo (M, \z) \times ([0,1]^k)^\z$ 
 are finite on every $Y_i$ and hence $(\mu, f^\z)$ is $0$-dimensional. Then, by Theorem \ref{theorem-light-map},
$\mdim (Y, \z) \leq \mdim (M, \z) \times ([0,1]^k)^\z \leq \mdim (M, \z) +\mdim ([0,1]^k)^\z =k  $.  Thus we have
that  $ \mdim (Y,\z) =n\mdim (X, \z) \leq k= [nd]+1$  for every positive integer $n$ and hence $\mdim (X,\z) \leq d$
 and the theorem follows. $\black$

\end{section}

\begin{section}{Proof of Theorem \ref{theorem-extension}}
\label{proof-theorem} Following \cite{levin-1} we use the following   notation.
Let  $(Y, \re)$ be a dynamical system,  $\cal A$ a collection of subsets of $Y$ and $\alpha, \beta \in \re$ positive numbers.
The collection  $\cal A$ is said to be {\bf $(\alpha, \beta)$-fine} if  $\mesh({\cal A}+r)< \alpha$ for every $r \in [0, \beta]\subset  \re$.

We also need
\begin{proposition}
{\rm (\cite{levin-1})}
\label{prop-family} 
Let $q>2$ be an integer. Then there is  a finite collection ${\cal E}$ of  disjoint closed intervals in $ [0, q) \subset \re$
such that $\cal E$  splits into the union ${\cal  E}={\cal E}_1 \cup \dots \cup {\cal E}_q$  of $q$ disjoint subcollections 
having the property that for every $ t \in \re$  the set  $t +\z\subset \re$  meets at least  $q-2$   subcollections ${\cal E}_i$
(a set meets a collection  if there is a point of the set that covered by the collection). Moreover,  we may assume
that $\mesh {\cal E}$ is as small as we wish. 
\end{proposition}
{\bf Proof of Theorem \ref{theorem-extension}.}
Let   $f=(f_1, \dots, f_k):X \lo [0, 1]^k$ be any  map and let $\epsilon>0$ and $\delta>0$
be such that under each $f_i$ the image of every subset of $X$ of $\diam <3\epsilon$ 
is  of $\diam < \delta$ in $[0,1]$. Our  goal is to approximate $f$ by a $\delta$-close map  $\psi$ 
such that  the  fibers of $(\mu, \psi^\z) : X \lo M \times ([0,1]^k)^\z$ contain at most
$\gamma=k/(k-d)$ points with pairwise distances larger than $3\epsilon$

Replacing  $d$ by a slightly larger real number keeping the value of $[k/(k-d)]$ unchanged we may assume  $\mdim X < d$.
 Take natural numbers $n$ and $q$ such that
$\mdim X < n/q < d$  and set $m=qk$. Then $n<qd<qk=m$ and $m/(m-n) \leq k/(k-d)$. 
By Theorem \ref{theorem-borel} we have $\mdim X \times_\z \re =\mdim X $. Then, assuming that $n$ is large enough, there is 
 a  Kolmogorov-Ostrand cover
  ${\cal F}$ of $X \times_\z \re$   such that  ${\cal F}$ is $(\epsilon, q)$-fine,
 ${\cal F}$ covers $X \times_\z \re$ at least $m-n$ times and
 ${\cal F}$ splits into ${\cal F}={\cal F}_1 \cup \dots \cup {\cal F}_{m}$ the union of finite  families of
 disjoint closed sets ${\cal F}_i$. 
 
 Recall that we consider $X \times_\z \re$  with the natural  action of $\re$ and identify  $X$  with a  subset of $X \times_\z \re$
 as described in \cite{levin-1}.
 Let  $\sigma \in \re$ be such that $0< \sigma< \epsilon$ and   for every disjoint elements $F', F'' $ of $\cal F$  and $t\in [0,q]$
 the sets $F'+t$ and $F''+t$
 are $\sigma$-distant. By this  we mean that 
 no subset of $X \times_\z \re$  of $\diam < \sigma$ meets both   $F'+t$ and $F''+t$.
 
Since $\mdim X    < n/q$  there is  an  integer $l>1$  and an open cover of $X$ 
such that $\ord {\cal U} < ln$ and $\cal U$ is $(\sigma,lq)$-fine. 
 
Let $W_M$ be any  non-empty open set in $M$ such that the closures of $W_M +z$ are disjoint for the integers
$-4lq\leq z \leq 4lq$ and
 let $\xi_M : M\lo \re$ be a Lindenstrauss level function determined by $W_M $.
Denote $\xi =\xi_M \circ \mu : X \lo \re$, $W=\mu^{-1}(W_M)$, $W_{\dg\dg} =W +\z \cap [-2lq,2lq]$, $X^{\dg\dg}=X \setminus W_{\dg\dg}$
and note that the closures of $W+z$ in $X$  are disjoint  for  the integers $-4lq\leq z \leq 4lq$  and $\xi$ is a   level  function such that
$\xi(x+z)=\xi(x) +z$ for every  $x \in X^{\dg\dg}$ and  an integer $-2lq\leq z\leq 2lq$.

Consider  a finite  collection ${\cal E}$ of  disjoint closed intervals in $ [0, q) \subset \re$ 
satisfying  the conclusions of Proposition \ref{prop-family}. Let $\tau>0$ be such
that any distinct  elements of ${\cal E}+\z$  are  $\tau$-distant in $\re$.
Take any open cover $\cal V$ of $\re$ with $\ord {\cal V}=2$ and $\mesh {\cal V} < \tau$.
Then the cover ${\cal U} \vee \xi^{-1}({\cal V})$  can be refined by an open cover ${\cal O}$ such that
$\ord {\cal O} <  ln+1$.  Refine $\cal O$ by 
a  Kolmogorov-Ostrand cover
  ${\cal F}^W$ of $X $   such that 
 ${\cal F}^W$ covers $X$ at least $l(m-n)$ times and
 ${\cal F}^W$ splits into ${\cal F}^W={\cal F}_1^W \cup \dots \cup {\cal F}_{lm}^W$ the union of finite  families of
 disjoint closed sets ${\cal F}^W_i$.  Clearly ${\cal F}^W$ is $(\sigma, lq)$-fine and $\mesh \xi({\cal F}^W) < \tau$.

We recall  the following notation from \cite{levin-1}.
Let  $\cal A$ be  a collection of subsets of $X \times_\z \re$, $\cal B$ a collection of intervals in $\re$. 
For $B\in \cal B$ and $z \in \z$ consider the collection ${\cal A}+B$ restricted to $\xi^{-1}(B+qz)$ and
denote by ${\cal A}\oplus_\xi B$ the union of such collections for all $z\in \z$.   Now  denote
by ${\cal A}\oplus_\xi {\cal B}$ the union
of the collections ${\cal A}\oplus_\xi  B$ for all $ B \in {\cal B}$. Note that ${\cal A}\oplus_\xi {\cal B}$ is
a collection of subsets of $X$.

 For $1\leq i \leq k$  define the collection  ${\cal D}_i$  of subsets of $X$ as the union of the collections
${\cal F}_{i}\oplus_\xi {\cal E}_1$, ${\cal F}_{i+k}\oplus_\xi {\cal E}_2$, \dots, ${\cal F}_{i+(q-1)k}\oplus_\xi {\cal E}_q$. 
Note that assuming that $\mesh \cal E$ is small enough we may also assume that ${\cal F}_i^+={\cal F}_i+[-\mesh{\cal E}, \mesh{\cal E}]$
is a collection of disjoint  $\sigma$-distant sets,
 the collection ${\cal F}^+={\cal F}+[-\mesh{\cal E}, \mesh{\cal E}]$ is 
$(\epsilon,q)$-fine   and, as a result,  we get that
 ${\cal D}_i $ is a collection of disjoint  closed sets of $X$ of $\diam < \epsilon$. Set 
 ${\cal D}$ to be the union of ${\cal D}_i$ for all $i$ and let
 ${\cal D}^{\dg}$  and ${\cal D}^{\dg}_i$  be the collections  $\cal D$ and  ${\cal D}_i$ respectively restricted to $X^{\dg}=X\setminus W_{\dg}$
 where $W_{\dg}=W+\z\cap [-lq,lq]$.

Now for every  $1\leq i \leq k$  define the collection  ${\cal D}^W_i$ of subsets of $X$ as the union of  the following  collections
for all integers $-2\leq z \leq 1$:
${\cal F}^W_i$  restricted to the closure of $W+zlq$, ${\cal F}^W_{i+k} +1$ restricted to the closure of $W+zlg+1$, $\dots$,
${\cal F}^W_{i+(lq-1)k} +lq$ restricted  to the closure of $W+zlq+lq-1=W+(z+1)ql-1$.
Set ${\cal D}^W$ to be the union of ${\cal D}^W_i$ for all $i$.
Clearly    ${\cal D}^W_i$ is a finite collection of disjoint closed sets and $\mesh {\cal D}^W< \sigma$.
Note that for every $D \in {\cal D}^W$ we have $\diam \xi(D\cap X^{\dg})< \tau$ because
$\mesh \xi({\cal  F}^W) < \tau$ and  $\xi(A+z)=\xi(A)+z$ for  every $A\subset X^\dg$ and $-lq\leq z \leq lq$.

Let us show that no element of ${\cal D}^W_i$  meets  disjoint elements of ${\cal D}^\dg_i$.
Let $D\in {\cal D}^W_i$ and let $D', D'' \in {\cal D}^{\dg}_i$ be disjoint. Recall that $\xi(D')$ and $\xi(D'')$ are contained
in the elements of ${\cal E}+\z$ and the elements of ${\cal E} +\z$ are pairwise disjoint.
Assume that there is $E \in \cal E$ such that
$\xi(D'), \xi(D'') \subset E+z$  for some $z\in \z$. Then $D'$ and $D''$ are $\sigma$-distant because
they are contained  in disjoint (and therefore $\sigma$-distant)  elements of ${\cal F}^+_i +t$ for some $t\in [0,q]$. 
 Since   $\diam D < \sigma$,
the set $D$  cannot meet both $D'$ and $D''$. Now assume that $\xi(D') \subset E'+z'$ and $\xi(D'') \subset E''+z''$
such that $E',E'' \in \cal E$, $z',z'' \in \z$ and $E'+z'$ and $E''+z''$ are disjoint. Hence $\xi(D')$ and $\xi(D'')$ are
$\tau$-distant in $\re$ and $D$ cannot meet both $D'$ and $D''$ because $\diam \xi(D\cap X^{\dg}) < \tau$.

Then there is a map $\psi=(\psi_1, \dots, \psi_k) : X\lo [0,1]^k$  so that  for  each $i$ the map  $\psi_i$ is $\delta$-close to $f_i$,
$\psi_i$ separates  the elements of ${\cal D}^W_i$ and $\psi_i$ also  separates the elements of  ${\cal D}_i^\dg$.  
Indeed,  first consider a  map $\psi_i$ which is $\delta$-close to $f_i$ and 
sends   the elements ${\cal D}_i^\dg$  to different singletons in $[0,1]$  and also sends the elements of ${\cal D}^W_i$
to (not necessarily different) singletons in $[0,1]$. Since no element of ${\cal D}^W_i$ meets
distinct elements of ${\cal D}_i^\dg$ we can replace $\psi_i$ by its approximation sending the elements of ${\cal D}^W_i$
to different singletons in $[0,1]$ preserving the property that $\psi_i$ separates the elements of ${\cal D}_i^\dg$.

We will show that the fibers of $(\mu, \psi^\z) : X \lo M\times [0,1]^\z$ contain 
at most  $\gamma$ points with pairwise distances  larger than $  \epsilon$. 
Aiming at  a contradiction assume that such  a  set  $\Gamma \subset  X$  exists with $|\Gamma |=[\gamma]+1$.
First  note  that    $\mu(\Gamma)$ is a singleton  in $M$ and therefore 
either $\Gamma \subset X^{\dg\dg}$ or $\Gamma \subset X \setminus X^{\dg\dg}=W_{\dg\dg}$. Also recall that $\xi$ factors through $\mu$
and therefore $\xi(\Gamma)$ is a singleton in $\re$ as well and denote $t_\Gamma =\xi(\Gamma) \in \re$.

Assume that $\Gamma \subset X^{\dg\dg}$.   Let $z$ be a non-negative integer    such that
$ zq\leq  t_{\Gamma}\leq (z+1)q$. 
Denote $t_*= t_{\Gamma} -[t_\Gamma]\in [0,1)$ and $q_*=[t_{\Gamma}]-zq\in [0,q)$ and note
that $0\leq t_*+q_* \leq q$.
 Consider  the set $\Gamma_* = \Gamma -q_*=\Gamma +zq+t_* \subset X$. 
Denote by $S$ the set of all the pairs   of integers  $(i,j)$  with $0 \leq i \leq q-1$ and $1\leq j \leq k$.
We  say that a point $x \in {\Gamma_*}$ is marked by a pair $(i,j) \in S$ if $x+i$ is covered by ${\cal D}_j$ 
and denote by $S_x \subset S$  the set of  the pairs by which $x$ is marked.
Let us compute the size of $S_x$.  Recall  that $x-t_*$ is covered by at least $m-n$ collections 
from the family ${\cal F}_1, \dots, {\cal F}_m$ and $t_*+\z$ meets   at least 
$q-2$ collections from ${\cal E}_1, \dots, {\cal E}_q$. 
Then  the point      $x $ is marked by at least $m-n-2k$ pairs $(i,j)$ of $S$ and hence  $|S_x| \geq m-n-2k$.

Indeed, for every ${\cal E}_p$ that meets  $t_* +\z$ pick up $i_p \in \z\cap [0,q)$ such that $t_*+i_p$ is covered by
${\cal E}_p$.  
 Note that  different $p$ define different $i_p$ and
for every $1\leq j\leq k$ such that ${\cal F}_{j+(p-1)k}$ covers the point  $x-t_*$ we have 
that the collection ${\cal D}_j$ covers $x+i_p$, and therefore $x$ is marked by the pair $(i_p,j)$ in $S$.
Thus if $t_*+\z$ meets all the collections ${\cal E}_1, \dots, {\cal E}_q$  
the number of   pairs $(i,j) \in S$  marking $x$ will be 
at least the number of times $x -t_*$ is covered 
by the collections ${\cal F}_1, \dots, {\cal F}_m$,  which is
at least $m-n$.
Each time  $t_*+\z$ misses a collection from ${\cal E}_1, \dots, {\cal E}_q$  reduces the above estimate by at most $k$. 
Since $t_*+\z$ can miss at most two collections from ${\cal E}_1, \dots, {\cal E}_q$
the point $x$ is marked by at least $m-n-2k$ pairs of $S$ and hence $|S_x|\geq m-n-2k$.

Since  $m-n\geq m(k-d)/k=q(k-d)$ we have 
$|S_x|\geq m-n-2k \geq q(k-d)-2k =q(k-d)\Delta$
where $\Delta =1-\frac{2k}{q(k-d)}$.
Note that we may take $q$ sufficiently large  and assume
that $\Delta < 1$ and 
 $\Delta > \gamma/|\Gamma|=\gamma/|\Gamma_*|$ and  get 
that  $\sum_{x \in \Gamma_*} |S_x | >|\Gamma_*| q(k-d)(\gamma/|\Gamma_*|)=qk$. Then, since $|S|=qk$,
there are two distinct points $x$ and $y$ in $\Gamma_*$  marked by the same  pair  $(i,j)$  of $S$ and hence
$x+i$ and $y+i$ are covered by ${\cal D}_j$.  Note that since  $x+q_*,y+q_* \in  \Gamma \subset X^{\dg\dg}$ we have  $x+i, y+i \in X^\dg$
and hence $x+i$ and $y+i$ are covered by ${\cal D}_j^\dg$. Recall that $\psi_j$ separates the elements of ${\cal D}_j^\dg$
and $x+i$ and $y+i$ are in the same fiber of $\psi_j$. Thus we get that $x+i$ and $y+i$ are contained in an element of ${\cal D}_j^\dg$
and hence in an element of ${\cal D}_j$.

Then $x-t_*$ and $y-t_*$ belong to the same element of ${\cal F}^+$ and hence the points
$x+q_*=(x-t_*)+(t_*+q_*)$ and $y+q_*=(y-t_*)+(t_*+q_*)$ of $\Gamma$
are $\epsilon$-close since ${\cal F}^+$ is $(\epsilon,q)$-fine and $0 \leq t_*+q_* \leq q$. 
Contradiction.

Now assume that $\Gamma \subset W_{\dg\dg}$.  We treat this case similarly to the previous one just adjusting 
the notations.  Since $\mu(\Gamma)$ is a singleton in $M$ we have
that $\Gamma \in W+zlq+p_*$ for some integers $z$ and $p_*$ such that  $-2\leq z \leq 1$ and $0\leq p_* \leq lq$.
Denote $\Gamma_*=\Gamma -p_*$. Let $S$ be the set of pairs $(i,j)$ with $0\leq i \leq lq-1$ and $1\leq j \leq k$
and say that a point $x \in \Gamma_*$ is marked by a pair $(i,j)$ of $S$ if $x+i$ is covered by ${\cal D}^W_j$.
Denote by $S_x$ the set of the pairs  of $S$  that mark $x$ and let us compute the size of $S_x$.
The point $x$ belongs to at least $l(m-n)$ families from  ${\cal F}^W_1, \dots , {\cal F}^W_{lm}$ and
for each such family ${\cal F}^W_{ik +j}$ with $0\leq i \leq lq-1$ and $1\leq j\leq k$ the point $x$
is marked by the pair $(i,j)$ of $S$. Thus we get that $|S_x| \geq l(m-n)\geq lq(k-d)$.
Since $|S|=klq$ and $|\Gamma_*|=|\Gamma|=[k/(k-d)]+1>k/(k-d)$ we get that  there are distinct points
$x$ and $y$ in $\Gamma_*$  marked by the same pair $(i,j)$ and hence
$x+i$ and $y+i$ are  covered by ${\cal D}^W_j$.  Recall that  $\psi_j$ separates  the elements of ${\cal D}_j^W$
and $x+i$ and $y+i$ are contained in the same fiber of $\psi_j$ and hence $x+i$ and $y+i$ are contained
in the same element of ${\cal D}^W_j$. Then $x$ and $y$ are contained in an element of ${\cal F}^W$ and
hence the points $x+p_*$ and $y+p_*$ of $\Gamma$ are $\epsilon$-close because
${\cal F}^W$ is $(\sigma, lq)$-fine, $0\leq p_*\leq lq$ and $\sigma < \epsilon$. Contradiction.

Thus 
  $\psi$ is the desired approximation of  $f$  the theorem follows by a standard Baire category  argument.
$\black$

\end{section}

 Department of Mathematics\\
Ben Gurion University of the Negev\\
P.O.B. 653\\
Be'er Sheva 84105, ISRAEL  \\
 mlevine@math.bgu.ac.il\\\\ 
  
\end{document}